\documentclass[10pt,english,reprint, superscriptaddress, secnumarabic, nofootinbib, nobibnotes]{revtex4-1}
\usepackage[T1]{fontenc}
\usepackage[latin9]{inputenc}
\usepackage[letterpaper]{geometry}
\usepackage{amsmath}
\usepackage{amsthm}

\makeatletter
\usepackage[all]{xy}

\usepackage{placeins}

\usepackage{lineno}

\def\frontmatter@abstractheading{}

\usepackage{babel}

\makeatother

\usepackage{babel}
\begin{document}

\title{Stochastic Unrelatedness, Couplings, and Contextuality}

\author{Ehtibar N.\ Dzhafarov}

\email{E-mail: ehtibar@purdue.edu }

\selectlanguage{english}%

\affiliation{Purdue University, USA }
\begin{abstract}
R. Duncan Luce once mentioned in a conversation that he did not consider
Kolmogorov's probability theory well-constructed because it treats
stochastic independence as a ``numerical accident,'' while it should
be treated as a fundamental relation, more basic than the assignment
of numerical probabilities. I argue here that stochastic independence
is indeed a ``numerical accident,'' a special form of stochastic
dependence between random variables (most broadly defined). The idea
that it is fundamental may owe its attractiveness to the confusion
of stochastic independence with stochastic unrelatedness, the situation
when two or more random variables have no joint distribution, ``have
nothing to do with each other.'' Kolmogorov's probability theory
cannot be consistently constructed without allowing for stochastic
unrelatedness, in fact making it a default situation: any two random
variables recorded under mutually incompatible conditions are stochastically
unrelated. However, stochastically unrelated random variables can
always be probabilistically coupled, i.e., imposed a joint distribution
upon, and this generally can be done in an infinity of ways, independent
coupling being merely one of them. The notions of stochastic unrelatedness
and all possible couplings play a central role in the foundation of
probability theory and, especially, in the theory of probabilistic
contextuality. 

\end{abstract}
\maketitle

\section{Introduction}

Almost 15 years ago R. Duncan Luce mentioned in a conversation that
the Kolmogorovian probability theory (KPT) was unsatisfactory because
it treated stochastic independence as a ``numerical accident'' rather
than a fundamental relation. If I roll a die {today }in Irvine,
California, Duncan said, and {on another day} you roll a die in
Lafayette, Indiana, the fact that the two outcomes are independent
cannot be established by checking the multiplication rule. On the
contrary, the applicability of the multiplication rule in this case
is justified by determining that the two dice are stochastically independent,
``have nothing to do with each other.''

This simple example (some may think too simple to be of great interest)
leads us to the very foundations of probability theory. Let us try
to understand it clearly by comparing it to another example. It is
about a situation when I repeatedly roll a single die, having defined
two random variables: 
\[
\begin{array}{l}
A=\left\{ \begin{array}{cc}
1 & \textnormal{if the outcome is even}\\
0 & \textnormal{otherwise}
\end{array}\right.,\\
\\
B=\left\{ \begin{array}{cc}
1 & \textnormal{if the outcome exceeds }3\\
0 & \textnormal{otherwise}
\end{array}\right..
\end{array}
\]
These two random variables co-occur in the most obvious empirical
meaning: the values of $A$ and $B$ are always observed together,
at every roll of the die. Another way of looking at it, the two random
variables co-occur because they are functions of one and the same
``background'' random variable $Z$, the outcome of rolling the
die. As a result, I can estimate from the observations the probabilities
$\Pr\left[A=1\textnormal{ and }B=1\right]$, $\Pr\left[A=1\right]$,
and $\Pr\left[B=1\right]$ (I will use $\Pr$ as a symbol for probability
throughout this paper): if the joint probability turns out to be the
product of the two marginal ones (statistical issues aside), the two
events are determined to be independent. I cannot simply make this
determination a priori, as it depends on what die I am rolling: if
it is a fair die, $A$ and $B$ are not independent, but if the distribution
of the outcomes is 
\[
\begin{array}{ccccccc}
value: & 1 & 2 & 3 & 4 & 5 & 6\\
pr.mass: & 0 & \frac{1}{4} & \frac{1}{4} & \frac{1}{4} & \frac{1}{4} & 0
\end{array},
\]
then $A$ and $B$ are independent.

The difference between this example and that of Duncan Luce's is not
in the number of the dice being rolled: my example would not change
too much if I roll two dice together, having marked them ``Left''
and ``Right,'' and define the random variables as 
\[
\begin{array}{l}
A=\left\{ \begin{array}{cc}
1 & \textnormal{if the Left outcome is even}\\
0 & \textnormal{otherwise}
\end{array}\right.,\\
\\
B=\left\{ \begin{array}{cc}
1 & \textnormal{if the Right outcome exceeds }3\\
0 & \textnormal{otherwise}
\end{array}\right..
\end{array}
\]
The realizations of $A$ and $B$ again come together, this time the
empirical meaning of the ``togetherness'' being ``in the same trial,''
or ``simultaneously.'' Again, one can also say that the two random
variables co-occur because they are functions of one and the same
``background'' random variable $Z$, only this time it is the pair
of values rather than a single one. And again, I can estimate from
the observations the probabilities $\Pr\left[A=1\textnormal{ and }B=1\right]$,
$\Pr\left[A=1\right]$, and $\Pr\left[B=1\right]$ and check their
adherence to the multiplication rule. Whether the two random variables
are stochastically independent is determined by the outcome of this
test: the dice may very well be rigged not to be independent.

In Duncan Luce's example the situation is very different: the outcomes
of rolling the two dice in two different places at two different times
have no empirically defined pairing. If I define my random variables
as 
\[
\begin{array}{l}
A=\left\{ \begin{array}{cc}
1 & {\textnormal{if on Tuesday in Irvine the outcome is even}}\\
0 & \textnormal{otherwise}
\end{array}\right.,\\
\\
B=\left\{ \begin{array}{cc}
1 & {\textnormal{if on Friday in Lafayette the outcome exceeds }3}\\
0 & \textnormal{otherwise}
\end{array}\right.,
\end{array}
\]
then I can estimate empirically the probabilities $\Pr\left[A=1\right]$,
and $\Pr\left[B=1\right]$ and find out, e.g., that they are (statistical
issues aside) 0.7 and 0.5, respectively. But I cannot estimate empirically
$\Pr\left[A=1\textnormal{ and }B=1\right]$: the two random variables
are not recorded in pairs. {The experiment involves no empirical
procedure by which one could find which value of $B$ should be paired
with which value of $A$. The two random variables therefore }do
not have an observable (estimable from frequencies) joint distribution,
they cannot be presented as functions of one and the same ``background''
random variable. What one can do, however, is to \emph{declare} the
two random variables stochastically independent, based on one's understanding
that they ``have nothing to do with each other.'' If one does so,
the validity of $\Pr\left[A=1\textnormal{ and }B=1\right]$ being
equal to the product of two individual probabilities is true \emph{by
construction}, requiring no empirical testing and allowing for no
empirical falsification.

This was Duncan Luce's point: while the KPT defines stochastic independence
through the multiplication rule, at least in some cases the determination
of independence precedes and justifies the applicability of the multiplication
rule. In Duncan Luce's opinion, this warranted treating stochastic
independence as a fundamental, ``qualitative'' relation preceding
assignment of numerical probabilities. This opinion is in accordance
with the general precepts of the representational theory of measurement.
Thus, the authors of the first volume of Foundations of Measurement
(Krantz et al., 1971) sympathetically refer to Zoltan Domotor 1969
dissertation in which he axiomatized probability theory treating stochastic
independence as a primitive relation. As far as I know, however, it
has not translated into a viable alternative to the KPT.

I accept Duncan Luce's example as posing a genuine foundational problem,
{but I disagree that this problem is about defining independence
by means other than the multiplication rule}. The position I advocate
below in this paper is as follows. 
\begin{enumerate}
\item Random variables that ``have nothing to do with each other'' are
defined on different {domains (}sample spaces). Rather than being
independent (which is a form of a joint distribution), they are \emph{stochastically
unrelated}, i.e., they possess no joint distribution. 
\item It is not that we do not know the ``true'' distribution, or that
in ``truth'' they are independent but we do not know how to justify
this. A joint distribution simply is not defined (until imposed by
us in one of multiple ways, discussed below). 
\item The KPT is consistent with the idea of multiple sample spaces and
in fact requires it for internal consistency: the idea of a single
sample space for all random variables imaginable is mathematically
untenable. 
\item Any given set of pairwise stochastically unrelated random variables
can always be \emph{coupled}, i.e., imposed a joint distribution on.
{This is equivalent to inventing a pairing scheme for their realizations,
}and this can be done in multiple ways, coupling them as independent
random variables being just one of them. 
\end{enumerate}

\section{On Random Variables, Unrelatedness, and Independence}

\subsection{Informal introduction}

Stochastic unrelatedness is easy to distinguish from stochastic independence:
the latter assumes the existence of a joint distribution, which means
that an empirical procedure exists by which each realization of one
random variables can be paired (coupled) with that of another. The
most familiar forms of empirical coupling are co-occurrence in the
same trial and co-relation to the same person. In the table below,
\begin{equation}
\begin{array}{ccccccc}
c: & 1 & 2 & 3 & 4 & 5 & \ldots\\
X: & x_{1} & x_{2} & x_{3} & x_{4} & x_{5} & \ldots\\
Y: & y_{1} & y_{2} & y_{3} & y_{4} & y_{5} & \ldots
\end{array},\label{eq: simple pairing}
\end{equation}
the indexing entity $c$ can be the number of a trial (as in repeatedly
rolling two marked dice together) or an ID of a person (as in relating
heights and weights, or weights before and after dieting). The random
variables $X$ and $Y$ here have a joint distribution: one can, e.g.,
estimate the probability with which $X$ falls within an event $E_{X}$
and (``simultaneously'') $Y$ falls within an event $E_{Y}$; and
if 
\begin{equation}
\Pr\left[X\in E_{X}\:\&\:Y\in E_{Y}\right]=\Pr\left[X\in E_{X}\right]\Pr\left[Y\in E_{Y}\right],\label{eq: simple independence}
\end{equation}
for any two such events $E_{X},E_{Y}$, then $X$ and $Y$ are considered
independent.

Suppose, however, that the information about $c$ in (\ref{eq: simple pairing})
{does not exist}, and {all one has} is some set of values for
$X$ and some set of values for $Y$. Clearly, now the ``togetherness''
of $X\in E_{X}$ and $Y\in E_{Y}$ is undefined. Although $\Pr\left[X\in E_{X}\right]$
and $\Pr\left[X\in E_{X}\right]$ have the same meaning as before,
$\Pr\left[X\in E_{X}\:\&\:Y\in E_{Y}\right]$ is undefined, and (\ref{eq: simple independence})
cannot be tested. This is what stochastic unrelatedness is: lack of
a joint distribution. A pair of stochastically unrelated random variables
are neither independent nor interdependent, these terms do not apply.

Think, e.g., of a list of weights {in some group of people} before
dieting ($X$) and a list of weights {in some other group of people
}after dieting ($Y$): {which value of $X$ should be paired with
which value of $Y$ to try to estimate $\Pr\left[X\in E_{X}\:\&\:Y\in E_{Y}\right]$?
Any pairing one can impose here will be as good as any other pairing,
and none of them is determined by the empirical procedures involved
(weighing people in the two groups). This simple example has counterparts
in all experiments where random variables are recorded under two or
more mutually exclusive conditions.}

A question arises: couldn't one nevertheless treat stochastically
unrelated $X$ and $Y$ \emph{as if} they were independent? The answer
is affirmative, but so is the answer to the question whether $X$
and $Y$ can be treated \emph{as if} they were not independent. Treating
the outcomes of rolling dice in {Irvine on Tuesday} ($X$) and in
{Lafayette on Friday} ($Y$) as if they were jointly distributed
means constructing another pair of random variables, $\left(\widetilde{X},\widetilde{Y}\right)$,
this time a jointly distributed one, such that $\widetilde{X}$ and
$\widetilde{Y}$ taken separately are distributed as $X$ and $Y$,
respectively. Such constructions form the subject of a special branch
of probability theory called the \emph{theory of coupling(s)} (Thorisson,
2000).

Let, e.g., both dice be fair. One can always construct $\left(\widetilde{X},\widetilde{Y}\right)$
by assigning probability mass $\frac{1}{36}$ to each of the 36 pairs.
This pair $\left(\widetilde{X},\widetilde{Y}\right)$ is the \emph{independent
coupling} of $X$ and $Y$. Its choice corresponds to pairing every
realization of $X$ with every realization of $Y$ {(or with uniformly
randomly chosen realization of $Y$)}. There is, however, no reason
to single out the independent coupling. One can also make $\widetilde{X}$
and $\widetilde{Y}$ perfectly correlated or perfectly anticorrelated
by assigning the probability masses as, respectively, 
\[
\mathrm{pr.mass}\left[\widetilde{X}=x\:\&\:\widetilde{Y}=y\right]=\begin{array}{ccc}
0 & if & x\not=y\\
\frac{1}{6} & if & x=y
\end{array}
\]
or 
\[
\mathrm{pr.mass}\left[\widetilde{X}=x\:\&\:\widetilde{Y}=y\right]=\begin{array}{ccc}
0 & if & x+y\not=7\\
\frac{1}{6} & if & x+y=7
\end{array},
\]
where $x,y\in\left\{ 1,\ldots,6\right\} $.

These couplings correspond to pairing each realization of $X$ with
only one specific realization of $Y$. If the dice producing $X$
and $Y$ have outcomes with different distributions, a perfectly correlated
or anticorrelated coupling will not be possible, while the independent
coupling will, as it is universally applicable. But the independent
coupling still will not be the only possible one (unless one of the
dice is rigged to roll a single outcome, in which case the only possible
coupling can be viewed as independent, perfectly correlated, or perfectly
anticorrelated).

{One may be tempted to think that the ``true'' pairing should involve
ordering the observations of $X$ and $Y$ chronologically and pairing
the outcomes with the same trial number. A brief reflection should
show, however, that this is an arbitrary choice: what theoretical
principles would compel one to pair the first realization of $X$
with the first realization of $Y$ (occurring, in Duncan Luce's example,
at another time in another place), rather than with the tenth one,
the last one, or one having the same quantile rank? Recall also that
the chronological sequences need not be defined to begin with: instead
of rolling a single die repeatedly one could roll a large number of
identical dice and count the events. }

Summarizing, a joint distribution for empirically observed $X$ and
$Y$ exists only if there is an empirical procedure for coupling their
realizations, such as relating them to one and the same value of $c$
in (\ref{eq: simple pairing}). Otherwise $X$ and $Y$ are stochastically
unrelated. When they are, one can impose on them a joint distribution
by creating a coupling $\left(\widetilde{X},\widetilde{Y}\right)$
for $X$ and $Y$ ``on paper.'' The individual distributions of
stochastically unrelated $X$ and $Y$ do impose some constraints
on possible joint distributions of $\left(\widetilde{X},\widetilde{Y}\right)$,
but, except in degenerate cases, do not determine it uniquely. The
independent coupling is not the only possible coupling of stochastically
unrelated random variables.

\subsection{\label{sub: Formalizing-the-naive}Formalizing the ``naive'' account
of random variables}

Random variables are defined by their distributions (say probability
masses associated with every possible roll of a die) and, in order
to distinguish different variables having the same distribution, by
their unique names (e.g., ``the outcome of the die rolled in Irvine
on Tuesday''). On a more general level, the distribution of a random
variable called $X$ is a probability space $\left(S_{X},\Sigma_{X},\mu_{X}\right)$,
with the standard meaning of the terms: $S_{X}$ is the set of possible
values for $X$, $\Sigma_{X}$ is a sigma-algebra of subsets of $S_{X}$,
and $\mu$ a probability measure.\footnote{I could have said ``distribution is determined by $\left(S_{X},\Sigma_{X},\mu_{X}\right)$,''
but it is simpler to say ``distribution is $\left(S_{X},\Sigma_{X},\mu_{X}\right)$,''
as we do not have an independent general definition of a distribution. } For each element $E_{X}$ of $\Sigma_{X}$ (an event) we define the
probability of $X$ ``falling in $E_{X}$'' or ``satisfying $E_{X}$''
as 
\begin{equation}
\Pr\left[X\in E_{X}\right]=\mu_{X}\left(E_{X}\right).
\end{equation}
Given another random variable, called $Y$ and distributed as $\left(S_{Y},\Sigma_{Y},\mu_{Y}\right)$,
we say it is jointly distributed with $X$ if there is a random variable
$Z=\left(X,Y\right)$ whose name is ``ordered pair of $X$ and $Y$''
and whose distribution is $\left(S_{X}\times S_{Y},\Sigma_{X}\otimes\Sigma_{Y},\nu\right)$,
subject to 
\begin{equation}
\begin{array}{c}
\nu\left(E_{X}\times S_{Y}\right)=\mu_{X}\left(E_{X}\right),\\
\\
\nu\left(S_{X}\times E_{Y}\right)=\mu_{Y}\left(E_{Y}\right),
\end{array}\label{eq: formal coupling of 2}
\end{equation}
for any events $E_{X}\in\Sigma_{X}$ and $E_{Y}\in\Sigma_{Y}$. The
meaning of $\Sigma_{X}\otimes\Sigma_{Y}$ is the smallest sigma-algebra
on $S_{X}\times S_{Y}$ that contains pairwise products of events
in $\Sigma_{X}$ and $\Sigma_{Y}$.

If such a $Z=\left(X,Y\right)$ exists (is defined among the random
variables one considers), then the joint distribution of $X$ and
$Y$ is unique. The existence of this random variable, however, is
not established by a mathematical derivation from the properties of
$X$ and $Y$, it is determined by the existence of an empirical procedure
in which the realizations of $X$ and $Y$ are observed ``together.''
If $Z=\left(X,Y\right)$ does not exist, one can always construct
a coupling $\widetilde{Z}=\left(\widetilde{X},\widetilde{Y}\right)$
whose distribution is $\left(S_{X}\times S_{Y},\Sigma_{X}\otimes\Sigma_{Y},\nu\right)$,
subject to (\ref{eq: formal coupling of 2}). The only difference
(but a critical one) is that the name of this $\widetilde{Z}$ is
not ``ordered pair of $X$ and $Y$'' but ``ordered pair of $\widetilde{X}$
{[}whose distribution is the same as that of $X${]} and $\widetilde{Y}$
{[}whose distribution is the same as that of $Y${]}.'' Such a $\widetilde{Z}$
can be freely introduced and does not change the $X$ and $Y$ being
coupled; in fact, $\left(\widetilde{X},\widetilde{Y}\right)$ is stochastically
unrelated to $X$ and to $Y$.

All of this can be easily generalized to an arbitrary set of random
variables (see, e.g., Dzhafarov and Kujala, in press).

\subsection{Random variables and joint distributions in KPT}

The formal account just given is not that of the traditional KPT.
The latter begins with the notion of a s\emph{ample space, }$\left(S,\Sigma,\mu\right)$,\footnote{{Note the terminological variance: ``sample space'' is more often
than not used in the literature to designate just the set $S$ rather
than the entire domain probability space $\left(S,\Sigma,\mu\right)$.
I find it more in line with the general meaning of the terms ``set''
and ``space'' in mathematics to refer to $S$ as a sample set (or
set of possible outcomes). This set is promoted into a space by endowing
it with a structure, which in this case is provided by the sigma algebra
$\Sigma$ and the measure $\mu$. }} and defines a random variables $X$ as a measurable mapping of this
space into a\emph{ }measurable space $\left(S_{X},\Sigma_{X}\right)$,
i.e., a function $X:S\rightarrow S'$ such that $X^{-1}\left(E_{X}\right)\in\Sigma$
for any $E_{X}\in\Sigma_{X}$.\footnote{Kolmogorov (1933/1956) only considered the case when $S_{X}$ is a
subset of reals and $\Sigma_{X}$ is the Borel sigma-algebra restricted
to this subset. In this paper I use the term ``random variable''
in the broad sense, with no restrictions imposed on $\left(S_{X},\Sigma_{X}\right)$.
Some authors prefer to use the term ``random element'' or ``random
entity'' to designate random variables in the broad sense. } The mapping induces on the codomain space $\left(S_{X},\Sigma_{X}\right)$
a probability measure $\mu_{X}$, by the rule 
\begin{equation}
\Pr\left[X\in E_{X}\right]=\mu_{X}\left(E_{X}\right)=\mu\left(X^{-1}\left(E_{X}\right)\right),
\end{equation}
for any $E_{X}\in\Sigma_{X}$. The resulting triple $\left(S_{X},\Sigma_{X},\mu_{X}\right)$
is called the \emph{distribution of} $X$. If another measurable mapping
$Y$ is defined on the same sample space, mapping it into a codomain
space $\left(S_{Y},\Sigma_{Y}\right)$ and resulting in the distribution
$\left(S_{Y},\Sigma_{Y},\mu_{Y}\right)$, then their joint distribution
$\left(S_{X}\times S_{Y},\Sigma_{X}\otimes\Sigma_{Y},\nu\right)$
is derived from the relation 
\begin{equation}
\Pr\left[X\in E_{X}\:\&\:Y\in E_{Y}\right]=\mu\left(X^{-1}\left(E_{X}\right)\cap Y^{-1}\left(E_{Y}\right)\right),
\end{equation}
for any $E_{X}\in\Sigma_{X}$ and $E_{Y}\in\Sigma_{Y}$. Note that
unlike in the ``naive'' approach above, the joint distribution of
$X$ and $Y$ here is mathematically derived from their individual
definitions as measurable functions on the same sample space.

Clearly, any two random variables defined on the same sample space
are jointly distributed. This may create a temptation to assume the
existence of a common sample space {and a joint distribution (even
if unknown to us) for }{\emph{any two}}{{} random variables.
In turn, this would mean the existence of a common sample space }for
\emph{all possible} random variables, so that random variables in
any set under consideration possess a joint distribution, and this
distribution is unique. Kolmogorov's (1933/1956) book may seem to
reinforce this view, as it does not explicitly speak of multiple sample
spaces. I disagree with this interpretation, even if not uncommon
(see, e.g., the overview of interpretations in Khrennikov, 2009c).
Kolmogorov's monograph ties the notion of a sample space\footnote{Kolmogorov's terminology is not the same as the modern terminology
(or variant thereof) I use in this paper; in particular, he does not
speak of a sample space but of a ``basic set with an algebra of subsets.''} to ``a complex of conditions which allows of any number of repetitions''
(Kolmogorov, 1933/1956, $\mathsection2$ of Chapter 1): this can be
interpreted as a position very close to if not the same as the one
argued for below. Whatever the correct interpretation, however, the
notion of a single sample space for all random variables is untenable
as it contradicts the common mathematical practices in dealing with
random variables. In Dzhafarov and Kujala (2014a) we presented the
following two arguments demonstrating this.

First of all, for any choice of a universal sample space $\left(S,\Sigma,\mu\right)$
all random variables $X$ defined on it will have the cardinality
of their set of possible values, defined as $S_{X}=X\left(S\right)$,
less than or equal to the cardinality of $S$. There is, however,
no justification, empirical or mathematical, for limiting the cardinality
of the set $S_{X}$ of all possible values for a random variable $X$.\footnote{Kolmogorov (1933/1956) did not have to deal with this issue, as the
random variables in this book are confined to real-valued ones.}

Second, even if we confine our attention to very simple random variables
with one and the same distribution, there is no justification, empirical
or mathematical, for limiting this set in any way. One can always
add a new random variable to any given set thereof. Thus, given any
set $\mathcal{N}$ of unit-normally distributed random variables,
one can introduce a unit-normally distributed $Y{}_{\mathcal{N}}$
such that its correlation with any $X\in\mathcal{N}$ is zero. If
there were a universal sample space $\left(S,\Sigma,\mu\right)$,
then there would be a definite set $\mathcal{N}^{*}$ of all possible
unit-normally distributed random variables. But this would mean that
our $Y_{\mathcal{N}^{*}}$ would have to belong to this set, which
is impossible, as $Y_{\mathcal{N}^{*}}$ cannot have zero correlation
with itself.

To further appreciate the untenability of a universal sample space,
observe that the identity mapping from this space into itself is a
random variable, $R$. The idea of a universal sample space therefore
is equivalent to the existence of a random variable $R$ of which
all imaginable random variables are functions. This does not add a
new formal argument against the idea, but it seems especially demonstrative:
what this mysterious ``super-variable'' $R$ could be?

\subsection{\label{sub:A-revised-KPT}A reinterpreted (or revised?) KPT}

All these considerations lead us to a different picture of the KPT,
in which there are different, stochastically unrelated random variables
$R,R',R'',\ldots$ , corresponding to different, mutually exclusive
conditions under which they are observed;\footnote{The notation $R,R',R'',\ldots$ is informal and should not be interpreted
as indicating a countable set.} and for each of these random variables one can define various functions
of it, 
\begin{equation}
\begin{array}{c}
X=f\left(R\right),\;Y=g\left(R\right),\ldots\\
X'=f'\left(R'\right),\;Y'=g'\left(R'\right),\ldots\\
X''=f''\left(R''\right),\;Y''=g''\left(R''\right),\ldots
\end{array}
\end{equation}
so that any two random variables that are functions of one and the
same member of the set $R,R',R'',\ldots$ possess a joint distribution,
while any functions of two different members of this set do not. To
preserve Kolmogorov's definition of a random variable, the $R,R',R'',\ldots$
can be thought of as identity functions on their separate sample spaces.

This picture is a step in the right direction, but it is still flawed
if we think of $R,R',R'',\ldots$ as some fixed set comprising all
pairwise stochastically unrelated random variables. The reason for
this is that if $R,R',R'',\ldots$ were a fixed set, with the corresponding
sample spaces (which, since they are identity functions, are simultaneously
their distributions) 
\begin{equation}
\left(S_{R},\Sigma_{R},\mu_{R}\right),\;\left(S_{R'},\Sigma_{R'},\mu_{R'}\right),\;\left(S_{R''},\Sigma_{R''},\mu_{R''}\right),\ldots\label{eq: spaces for R,R',R''}
\end{equation}
then one could form a single random variables $R^{*}$ of which $R,R',R'',\ldots$
(hence also all other random variables imaginable) were functions.
This ``super-variable'' $R^{*}$ would have the set and sigma-algebra
that are products of, respectively, sets and sigma-algebras in (\ref{eq: spaces for R,R',R''}),
and it would have a probability measure $\nu$ from which $\mu_{R},\mu_{R'},\mu_{R''},\ldots$
are computed as marginals, e.g., 
\[
\mu_{R}\left(E_{R}\right)=\nu\left(E_{R}\times S_{R'}\times S_{R''}\times\ldots\right),
\]
for any $E_{R}\in\Sigma_{R}$. We have seen already that the idea
of such a ``super-variable'' is untenable.

{A logically consistent way out of this difficulty} is to consider
$R,R',R'',\ldots$ as a \emph{class with uncertain and/or flexible
membership}.\footnote{{A different approach is presented in Dzhafarov and Kujala (2015b),
where we formally define, by means of a quasi-constructive procedure,
the set of all random variables considered ``existing'' in a given
study. }} Indeed, it should be clear from the previous discussion that random
variables can be freely introduced, so, e.g., there is no fixed set
of random variables with any given distribution. Some random variables
we observe have an empirically defined coupling scheme, and then they
are jointly distributed. Other sets of random variables we observe
are observed under different conditions each and do not have an empirical
coupling. Then they can be modeled as stochastically unrelated random
variables. However, we then can create ``copies'' of these random
variables and couple them ``on paper'' in a multitude of ways. This
seems to be a consistent view of random variables. KPT is by no means
dismissed in this view, because any distribution $\left(S_{X},\Sigma_{X},\mu_{X}\right)$
for a random variable $X$ is a probability space subject to Kolmogorov's
axioms: 
\begin{enumerate}
\item $\mu_{X}$ is a function $\Sigma_{X}\rightarrow\left[0,\infty\right)$; 
\item $\mu_{X}$$\left(S_{X}\right)=1$; 
\item $\mu_{X}$$\left(\bigcup_{i=1}^{\infty}E_{X}^{\left(i\right)}\right)=\sum_{i=1}^{\infty}\mu_{X}\left(E_{X}^{\left(i\right)}\right)$
for any sequence of pairwise disjoint $E_{X}^{\left(1\right)},E_{X}^{\left(2\right)},\ldots$
in $\Sigma_{X}$. 
\end{enumerate}
Moreover, insofar as one focuses on a given set of jointly distributed
random variables, all of them can be presented as measurable functions
on a single sample space (or functions of a single random variable).

\subsection{Radical contextualism}

Is there a unique way of determining which random variables are and
which are not stochastically interrelated? A general answer to this
question is negative: the definition of a jointly distributed set
of random variables involves an empirical procedure of coupling their
realizations, ``observing them together.'' The meaning of such an
empirical procedure may be different for different situation and different
observers. From the mathematical point of view, however, the question
is about a language that makes the fundamental distinctions between
stochastically related and stochastically unrelated random variables.
Such a language is proposed in Dzhafarov and Kujala (2014b), Dzhafarov,
Kujala, and Larsson (2015), Kujala and Dzhafarov (2015), and Kujala,
Dzhafarov, and Larsson (2015). For an overview, see Dzhafarov and
Kujala (2015a) and Dzhafarov, Kujala, and Cervantes (2016).

It is postulated that every random variable's identity is determined
by two types of variables, referred to as \emph{objects} (also properties,
entities, contents, etc.) and \emph{contexts }(also conditions, environment,
etc.). Intuitively, the random variables are treated as ``measurements,''
the objects answer the question ``what is measured?'' whereas the
contexts answer the question ``how is it measured?''

Let $Q$ be a set of objects and $C$ a set of contexts considered
in a given study. The mentioning of ``a given study'' is essential:
in a different study one could choose a different set of objects to
measure and a different set of contexts in which to measure them.
The measurement $R$ of an object $q\in Q$ in a context $c\in C$
is denoted by $R_{q}^{c}$.

The meaning of a context is that it provides an empirical coupling
for the measurements within this context: the random variables $R_{q}^{c}$
with different $q$ measured within the same context $c$ are ``measured
together,'' i.e., they possess a joint distribution. Denoting by
$Q_{c}$ the subset of objects in $Q$ that are measured in context
$c$, 
\begin{equation}
R^{c}=\left\{ R_{q}^{c}\right\} _{q\in Q_{c}}\label{eq: R^c}
\end{equation}
is a random variable (which implies that it has a distribution, and
this distribution is a joint distribution of its components). On the
other hand, any two random variables $R_{q}^{c}$ and $R_{q'}^{c'}$
with $c\not=c'$ are stochastically unrelated, whether $q$ and $q'$
are distinct objects or not. It follows that any random variables
$R^{c}$ and $R^{c'}$ defined as in (\ref{eq: R^c}) with $c\not=c'$
are stochastically unrelated.

The idea of such contextual notation and (mutatis mutandis) understanding
of stochastic unrelatedness have precursors and analogues in the quantum-physical
literature: see Khrennikov (2005, 2008, 2009a-c), Simon, Brukner,
and Zeilinger (2001), Larsson (2002), Svozil (2012), and Winter (2014).
Khrennikov (2009c) points out that the contextual understanding of
random variables is intrinsic features of von Mises's ``ensemble
approach'' to probabilities: the identity of an ``ensemble'' of
observations corresponds to context in which these observations are
made.\footnote{Khrennikov thinks that in this respect von Mises's approach is radically
different from Kolmogorov's, an opinion one can disagree with if the
KPT is not confined to a single probability space.}

Three aspects of our theory, however, set it aside from this literature: 
\begin{enumerate}
\item Contextual labeling is universal, and no two random variables recorded
in different contexts have a joint distribution. 
\item Pairwise stochastically unrelated random variables $\left\{ R^{c}\right\} _{c\in C}$
(each of which is a set of jointly distributed random variables) can
be coupled at will, with no coupling being privileged. 
\item The random variables $\left\{ R^{c}\right\} _{c\in C}$ can be characterized
by whether it is possible or impossible to couple them in a particular
way (e.g., by a maximally connected coupling, as discussed in Section
\ref{sub: An-example-of}). 
\end{enumerate}
Below I will give an example of how these principles work in solving
the problem of selective influences in psychology, as well as its
generalized version, the problem of contextuality, primarily studied
in quantum mechanics. First, however, I have to address some obvious
objections to the radical contextualism.

\subsection{Possible objections}

The first objection is that it is impossible to take into account
all conditions in the world, and without knowing them one would not
know if one deals with stochastically related or unrelated random
variables. The response to this objection lies in the qualification
``in a given study'' I made when I introduced object sets $Q$ and
context sets $C$. The identification of random variables by what
they measure and by how they measure it depends on what other variables
in the world one records and relates to realizations of the random
variables in question.

To give an example, let there be a very large group of husband-and-wife
couples; to each of the husbands Alice poses one of two different
Yes/No questions, $a_{1}$ or $a_{2}$; to each of the wives Bob poses
one of two different Yes/No questions, $b_{1}$ or $b_{2}$ (that
may be the same as or different from $a_{1},a_{2}$). Alice decides
(this is not a matter of truth or falsity but one of convention) to
consider the responses to $a_{1}$ and $a_{2}$ in the group of husbands
as realizations of random variables $R_{a_{1}}$ and $R_{a_{2}}$,
respectively; and Bob defines $R_{b_{1}}$ and $R_{b_{2}}$ for the
group of wives analogously. This labeling indicates that Alice treats
$a_{1}$ and $a_{2}$ as objects being measured (by responses to these
questions), and so does Bob for $b_{1}$ and $b_{2}$.\footnote{This example is formally equivalent to the EPR/Bohm experiment in
quantum physics (see Section \ref{sub: An-example-of}).}

Let us ask now: what are the contexts in which Alice records her $R_{a_{1}}$
and $R_{a_{2}}$? By the rules of the survey each person answers a
single question, so asking $a_{1}$ excludes asking $a_{2}$ and vice
versa. In other words, the conditions under which one records answers
to $a_{1}$ and $a_{2}$ are incompatible. Formally, this means that
$a_{1}$ is measured in the context $a_{1}$ while $a_{2}$ is measured
in the context $a_{2}$: Alice has therefore stochastically unrelated
random variables $R_{a_{1}}^{a_{1}}$ and $R_{a_{2}}^{a_{2}}$ . Analogously,
Bob has stochastically unrelated $R_{b_{1}}^{b_{1}}$ and $R_{b_{2}}^{b_{2}}$.
Their stochastic unrelatedness is quite obvious: why should any response
by a person to $a_{1}$ (or $b_{1}$) be paired with any response
by another person to $a_{2}$ (respectively, $b_{2}$)?

By the same argument, either of Alice's measurements is stochastically
unrelated to either of Bob's: the four measurements 
\begin{equation}
R_{a_{1}}^{a_{1}},R_{a_{2}}^{a_{2}},R_{b_{1}}^{b_{1}},R_{b_{2}}^{b_{2}}\label{eq: 4 separate contexts}
\end{equation}
are made in four different contexts. It is clear, however, that Alice
and Bob could try to form joint distributions of their measurements
using some empirical coupling procedure, e.g., the pairing of the
measurements by the marital relation: that is, pairing a husband's
response to $a_{i}$ with his wife's response to $b_{j}$, for each
of the four combinations of $i=1,2$ and $j=1,2$. To do this means
to form new contexts, $\left(a_{1},b_{1}\right)$, $\left(a_{1},b_{2}\right)$,
$\left(a_{2},b_{1}\right)$, and $\left(a_{2},b_{2}\right)$, and
to re-label Alice's random variables as 
\begin{equation}
R_{a_{1}}^{\left(a_{1},b_{1}\right)},R_{a_{1}}^{\left(a_{1},b_{2}\right)},R_{a_{2}}^{\left(a_{2},b_{1}\right)},R_{a_{2}}^{\left(a_{2},b_{2}\right)}
\end{equation}
while Bob's random variables become 
\begin{equation}
R_{b_{1}}^{\left(a_{1},b_{1}\right)},R_{b_{1}}^{\left(a_{2},b_{1}\right)},R_{b_{2}}^{\left(a_{1},b_{2}\right)},R_{b_{2}}^{\left(a_{2},b_{2}\right)}.
\end{equation}
The previous four pairwise stochastically unrelated variables (\ref{eq: 4 separate contexts})
are replaced now with the four pairwise stochastically unrelated variables
\begin{equation}
R^{\left(a_{i},b_{j}\right)}=\left(R_{a_{i}}^{\left(a_{i},b_{j}\right)},R_{b_{j}}^{\left(a_{i},b_{j}\right)}\right),\;i,j\in\left\{ 1,2\right\} .\label{eq: 4 bunches}
\end{equation}
There is no justification for saying that either of these representations,
(\ref{eq: 4 separate contexts}) or (\ref{eq: 4 bunches}), is more
``correct'' than another: Alice who does not know whose husband
her respondent is and Alice who knows this deal with different sets
of random variables.

Another, related objection is that radical contextualism should lead
to considering every realization of a random variable as being stochastically
unrelated to every other realization. In the previous example, if
Alice records the identities of the people she is posing the question
$a_{1}$ to, then in place of a single $R_{a_{1}}^{a_{1}}$ she creates
random variables $R_{a_{1}}^{John}$, $R_{a_{1}}^{Peter}$, etc.,
each with a single realization. In a typical behavioral experiment,
John, Peter, etc. can be replaced with ``trial 1,'' ``trial 2,''
etc. As the contexts (upper indexes) differ, the variables are stochastically
unrelated. Isn't this a problem? In particular, does not stochastic
unrelatedness of the realizations of a random variable clash with
the standard statistical practice of viewing them as independent identically
distributed variables?

The answer to these questions is essentially the same as to the previous
objection. Alice does not have to record the identity of the people
she queries, and if she does not, then $R_{a_{1}}^{a_{1}}$ is the
random variable she forms. If she does the recording, then she creates
new contexts, and it is indeed true then that $R_{a_{1}}^{John}$,
$R_{a_{1}}^{Peter}$, etc, are pairwise stochastically unrelated.
They can, however, be coupled, and as always when coupling is not
based on an empirical procedure, this can be done in a multitude of
ways. One possible coupling is the independent coupling, the other
is the identity coupling, and one could create an infinity of other
couplings.

This may be difficult to understand. Suppose we know that John responded
Yes, Peter responded No, Paul responded No, etc. --- then how could
one speak of the identity coupling? Or, if we know that the response
in trial $n+1$ repeats the response in trial $n$ with probability
0.7 --- how can one speak then of the independent coupling? To answer
these questions one should recall that by coupling pairwise stochastically
unrelated $R^{\left(1\right)},R^{\left(2\right)},R^{\left(3\right)},\ldots$
one does not mysteriously transform them into jointly distributed
random variable. Instead one creates a new sequence $\left(\widetilde{R}^{\left(1\right)},\widetilde{R}^{\left(2\right)},\widetilde{R}^{\left(3\right)},\ldots\right)$,
in which each $\widetilde{R}^{\left(i\right)}$ has the same distribution
as $R^{\left(i\right)}$, and all these components have a joint distribution.
The table below demonstrates the logical structure of the identity
coupling:

\begin{equation}
\begin{array}{ccccc}
 & R^{\left(1\right)} & R^{\left(2\right)} & R^{\left(3\right)} & \ldots\\
\widetilde{R}^{\left(1\right)} & \boxed{\overset{}{\:r_{1}}\:} & r_{2} & r_{3} & \ldots\\
\widetilde{R}^{\left(2\right)} & r_{1} & \boxed{\overset{}{\:r_{2}}\:} & r_{3} & \ldots\\
\widetilde{R}^{\left(3\right)} & r_{1} & r_{2} & \boxed{\overset{}{\:r_{3}}\:} & \ldots\\
\vdots & \vdots & \vdots & \vdots & \ddots
\end{array}
\end{equation}
The boxed values are the ones factually observed, the rest of the
values in each column are those attained by the corresponding components
of the identity coupling. As we see, this has nothing to do with the
observed values being or not being equal to each other.

The next table demonstrates the logical structure of the independent
coupling: 
\begin{equation}
\begin{array}{ccccc}
 & R^{\left(1\right)} & R^{\left(2\right)} & R^{\left(3\right)} & \ldots\\
\widetilde{R}^{\left(1\right)} & \boxed{\overset{}{\:r_{1}}\:} & r'_{2} & r'_{3} & \ldots\\
\widetilde{R}^{\left(2\right)} & r'_{1} & \boxed{\overset{}{\:r_{2}}\:} & r''_{3} & \ldots\\
\widetilde{R}^{\left(3\right)} & r''_{1} & r''_{2} & \boxed{\overset{}{\:r_{3}}\:} & \ldots\\
\vdots & \vdots & \vdots & \vdots & \ddots
\end{array}
\end{equation}
The boxed values, again, are those factually observed, and the primed
values in the $i$th column are sampled from a coupling $\left(\widetilde{R}^{\left(1\right)},\widetilde{R}^{\left(2\right)},\widetilde{R}^{\left(3\right)},\ldots\right)$
with stochastically independent components and $\widetilde{R}^{\left(i\right)}=r_{i}$.
Again, this has nothing to do with the observed values forming or
not forming a sequence with certain statistical properties.

Focusing on the statistical properties of the observed (``boxed'')
values means, formally, that the observations in different trials
(or responses from different persons) are treated as objects rather
than contexts, all these objects being measured in a single context
and therefore jointly distributed: 
\begin{equation}
R_{1}^{\left(1,2,3,\ldots\right)},R_{2}^{\left(1,2,3,\ldots\right)},R_{3}^{\left(1,2,3,\ldots\right)},\ldots
\end{equation}

A third objection one can raise against the radically contextual reinterpretation
(or revision) of the KPT is that the notions of an ``object'' and
a ``context'' are not mathematically defined: they are primitives
of the language proposed. How can one know what objects and what contexts
to invoke in a specific situation? The response to this is that it
is indeed not a mathematical issue. Mathematical analysis begins once
one has specified a set $Q$ of objects and a set $C$ of contexts,
and there is no single correct way of doing it.

Consider, e.g., the situation when two questions are asked in one
of two orders, $a\rightarrow b$ or $b\rightarrow a$. One can take
$a$ to be the same object measured in two different contexts, and
similarly for $b$, forming thereby four random variables (responses
to the questions) 
\begin{equation}
R_{a}^{a\rightarrow b},R_{b}^{a\rightarrow b},R_{a}^{b\rightarrow a},R_{b}^{b\rightarrow a}.
\end{equation}
By our rules, they are grouped into two stochastically unrelated random
variables 
\begin{equation}
R^{a\rightarrow b}=\left(R_{a}^{a\rightarrow b},R_{b}^{a\rightarrow b}\right)\textnormal{ and }R^{b\rightarrow a}=\left(R_{a}^{b\rightarrow a},R_{b}^{b\rightarrow a}\right).\label{eq: order effect normal}
\end{equation}
This view of the situation leads to an interesting contextual analysis
(Dzhafarov, Zhang, \& Kujala, 2015).

It is, however, possible to deny that ``the same question $a$''
means ``the same object $a$'' in the two contexts: one can maintain
instead that $a$ asked first is simply a different object from $a$
asked second; and similarly for $b$. In this view we have four objects,
$a_{1},a_{2},b_{1},b_{2}$ (where index indicates whether the question
is asked first or second), measured in two contexts, $a_{1}\rightarrow b_{2}$
and $b_{1}\rightarrow a_{2}$. One ends up with two stochastically
unrelated random variables 
\begin{equation}
\begin{array}{c}
R^{a_{1}\rightarrow b_{2}}=\left(R_{a_{1}}^{a_{1}\rightarrow b_{2}},R_{b_{2}}^{a_{1}\rightarrow b_{2}}\right)\\
\textnormal{and}\\
R^{b_{1}\rightarrow a_{2}}=\left(R_{a_{2}}^{b_{1}\rightarrow a_{2}},R_{b_{1}}^{b_{1}\rightarrow a_{2}}\right).
\end{array}\label{eq: order effect alternative}
\end{equation}
This representation allows for no nontrivial contextual analysis (see
below), as the stochastically unrelated random variables have no objects
in common. It is, nevertheless, as legitimate a representation as
the previous one. A psychologist will most probably choose (\ref{eq: order effect normal})
over (\ref{eq: order effect alternative}) (Wang \& Busemeyer, 2013;
Wang et al., 2014), but it is not mathematics that dictates this choice.

\subsection{\label{sub: An-example-of}An example of contextual analysis}

The problem of selective influences was introduced to psychology by
Sternberg (1969) and developed through a series of publications (Schweickert
\& Townsend, 1989; Townsend, 1984, 1990; Townsend \& Schweickert,
1989; Roberts \& Sternberg, 1993; Townsend \& Nozawa, 1995; Schweickert,
Giorgini, \& Dzhafarov, 2000; Dzhafarov 2003; Dzhafarov, Schweickert,
\& Sung, 2005; Kujala \& Dzhafarov, 2008; Dzhafarov \& Kujala, 2010).
Later, a link has been established between this problem and the quantum-mechanical
analysis of entanglement (Dzhafarov \& Kujala, 2012a-b, 2013, 2014c)
and, more generally, probabilistic contextuality (Dzhafarov \& Kujala,
2014a-b, 2015a-b; Dzhafarov, Kujala, \& Larsson, 2015; Kujala, Dzhafarov,
\& Larsson, 2015).

I will formulate the problem using the contextual language introduced
above. Let there be a system acted upon by two inputs, $\alpha$ and
$\beta$, and reacting by two simultaneous distinct responses, $R_{\alpha}$
and $R_{\beta}$ (or distinct aspects of the same response, such as
response time and response accuracy). The indexation here reflects
the belief (or hypothesis) that $R_{\alpha}$ is ``primarily'' influenced
by $\alpha$ and $R_{\beta}$ by $\beta$. One can also say that $R_{\alpha}$
measures $\alpha$ and $R_{\beta}$ measures $\beta$. The question
is whether $R_{\alpha}$ is also influenced by $\beta$ and/or $R_{\beta}$
is also influenced by $\alpha$. Let us simplify the problem by assuming
that $\alpha\in\left\{ 1,2\right\} $ and $\beta\in\left\{ 1,2\right\} $,
and they vary in a completely crossed factorial design, $\left\{ 1,2\right\} \times\left\{ 1,2\right\} $.
Each of the treatments $\left(\alpha,\beta\right)=\left(i,j\right)$
should be considered a context, wherefrom the responses of the system
must be labeled 
\begin{equation}
R^{\left(\alpha=i,\beta=j\right)}=\left(R_{\alpha=i}^{\left(\alpha=i,\beta=j\right)},R_{\beta=j}^{\left(\alpha=i,\beta=j\right)}\right),\;i,j\in\left\{ 1,2\right\} .
\end{equation}
To remind the interpretation, $R_{\alpha=i}^{\left(\alpha=i,\beta=j\right)}$
measures the object $\alpha=i$ in the context $\left(\alpha=i,\beta=j\right)$;
$R_{\beta=j}^{\left(\alpha=i,\beta=j\right)}$ measures the object
$\beta=j$ in the same context; being in the same context, these two
measurements form a random variable $R^{\left(\alpha=i,\beta=j\right)}$
(whose components possess a joint distribution); however, the four
random variables $R^{\left(\alpha=i,\beta=j\right)}$ are pairwise
stochastically unrelated. To lighten the notation, let us put 
\begin{equation}
R_{\alpha=i}^{\left(\alpha=i,\beta=j\right)}=A_{i}^{ij},\;R_{\beta=j}^{\left(\alpha=i,\beta=j\right)}=B_{j}^{ij}.
\end{equation}

According to the definition of selective influences given in Dzhafarov
(2003) and elaborated in Dzhafarov and Kujala (2010), one says that
$A_{i}^{ij}$ is not influenced by $\beta$ and $B_{j}^{ij}$ is not
influenced by $\alpha$ (for all $i,j$) if one can find a coupling
\begin{equation}
\left(\widetilde{A}_{1}^{11},\widetilde{B}_{1}^{11},\widetilde{A}_{1}^{12},\widetilde{B}_{2}^{12},\widetilde{A}_{2}^{21},\widetilde{B}_{1}^{21},\widetilde{A}_{2}^{22},\widetilde{B}_{2}^{22}\right)\label{eq: coupling 2x2}
\end{equation}
in which the equalities 
\begin{equation}
\begin{array}{ccc}
\widetilde{A}_{1}^{11}=\widetilde{A}_{1}^{12}, &  & \widetilde{A}_{2}^{21}=\widetilde{A}_{2}^{22},\\
\\
\widetilde{B}_{1}^{11}=\widetilde{B}_{1}^{21}, &  & \widetilde{B}_{2}^{12}=\widetilde{B}_{2}^{22}
\end{array}\label{eq: equalities}
\end{equation}
hold with probability 1. Put differently, the random variables $\widetilde{A}_{1}^{11}$
and $\widetilde{A}_{1}^{12}$ in the joint distribution of (\ref{eq: coupling 2x2})
always attain one and the same value, even though the value of $\beta$
changes; and analogously for the remaining three equalities. Note
that, by the definition of a coupling, 
\begin{equation}
\widetilde{A}_{i}^{ij}\stackrel{d}{=}A_{i}^{ij}\textnormal{ and }\widetilde{B}_{j}^{ij}\stackrel{d}{=}B_{j}^{ij},\;i,j\in\left\{ 1,2\right\} ,
\end{equation}
where $\stackrel{d}{=}$ means ``has the same distribution as.''

If all $A$ and $B$ responses of the system have a finite number
of possible values, this situation generalizes Bohm's version of the
Einstein-Podolsky-Rosen (EPR) paradigm (Bohm \& Aharonov, 1957; Bell,
1964). Of course, the distributions of the $A,B$ need not be generally
in compliance with the quantum rules for entangled particles, but
the existence or nonexistence of a coupling with the stipulated properties
should be determinable for any observed $A$ and $B$. Let us assume
for simplicity that both $A$ and $B$ responses of the system are
binary, and let us denote their values $+1$ and $-1$. In this special
case the necessary and sufficient conditions for the selectiveness
of influences are given by 
\begin{equation}
\begin{array}{ccc}
A_{i}^{i1}\stackrel{d}{=}A_{i}^{i2} & for & i=1,2\\
B_{j}^{1j}\stackrel{d}{=}B_{j}^{2j} & for & j=1,2
\end{array}\label{eq: marginal}
\end{equation}
and 
\begin{equation}
\max_{k,l\in\left\{ 1,2\right\} }\left|\sum_{i,j\in\left\{ 1,2\right\} }\mathsf{E}\left[A_{i}^{ij}B_{j}^{ij}\right]-2\mathsf{E}\left[A_{k}^{kl}B_{l}^{kl}\right]\right|\leq2,\label{eq: CHSH}
\end{equation}
where $\mathsf{E}$ stands for expected value. This (in an algebraically
different form) was first proved by Fine (1982), although (\ref{eq: marginal})
in his work is implied by the notation rather than stated explicitly.
The distributional equalities (\ref{eq: marginal}) describe the condition
known as \emph{marginal selectivity}: the distribution of $A_{i}^{ij}$
does not change with the value $j$ of $\beta$, and the distribution
of $B_{j}^{ij}$ does not change with the value $i$ of $\alpha$.
The numerical inequality (\ref{eq: CHSH}) is known as the \emph{CHSH
inequality} (after the authors of Clauser et al., 1969). In quantum
mechanics, violations of this inequality when the marginal selectivity
(\ref{eq: marginal}) holds is described by saying that the system
is \emph{contextual} (see, e.g., Kurzynski, Ramanathan, \& Kaszlikowski,
2012).

If the marginal selectivity (\ref{eq: marginal}) is violated, the
CHSH inequality (\ref{eq: CHSH}) cannot be derived, and it makes
no difference whether it is satisfied or not. Moreover, if marginal
selectivity is violated, it seems unnecessary to look at anything
else: clearly then $A$ is not selectively influenced by $\alpha$
alone, and/or $B$ is not selectively influenced by $\beta$ alone.
As it turns out, however, one may still be interested in the question:
is the influence of $\beta$ upon $A$ and/or of $\alpha$ upon $B$
entirely described by the violations of marginal selectivity? Indeed,
since the CHSH inequality (\ref{eq: CHSH}) may very well be violated
when the marginal selectivity (\ref{eq: marginal}) holds, and since
we then conclude that selectiveness of influences is violated too,
we have to admit that the ``wrong'' influences (from $\beta$ to
$A$ and/or from $\alpha$ to $B$) can be indirect, without manifesting
themselves in changed marginal distributions. This leads us to a generalized
notion of contextuality (Dzhafarov, Kujala, \& Larsson, 2015; Kujala,
Dzhafarov, \& Larsson, 2015; Dzhafarov \& Kujala, 2015a-b; Dzhafarov,
Zhang, Kujala, 2015; Dzhafarov, Kujala, \& Cervantes, 2016).

When applied to our example with two binary inputs $\alpha,\beta$
and two binary random outputs $A,B$, the definition is as follows.
A system 
\[
\left(A_{1}^{11},B_{1}^{11}\right),\left(A_{1}^{12},B_{2}^{12}\right),\left(A_{2}^{21},B_{1}^{21}\right),\left(A_{2}^{22},B_{2}^{22}\right)
\]
is noncontextual if it has a \emph{maximally connected coupling}.
The latter is defined as a coupling (\ref{eq: coupling 2x2}) in which
each of the equalities (\ref{eq: equalities}) holds with the maximal
possible probability that is allowed by the individual distributions
of the random variables. To explain, if $A_{1}^{11}\stackrel{d}{=}A_{1}^{12}$,
then the maximal possible value for $\Pr\left[\widetilde{A}_{1}^{11}=\widetilde{A}_{1}^{12}\right]$
is $1$. Applying this to all other equalities in (\ref{eq: equalities}),
we get the previous definition. If, however, $A_{1}^{11}$ and $A_{1}^{12}$
have different distributions, then the maximal possible value for
$\Pr\left[\widetilde{A}_{1}^{11}=\widetilde{A}_{1}^{12}\right]$ is
\begin{equation}
\begin{array}{l}
\min\left\{ \Pr\left[A_{1}^{11}=1\right],\Pr\left[A_{1}^{12}=1\right]\right\} \\
\\
+\min\left\{ \Pr\left[A_{1}^{11}=-1\right],\Pr\left[A_{1}^{12}=-1\right]\right\} \\
\\
=1-\left|\Pr\left[A_{1}^{11}=1\right]-\Pr\left[A_{1}^{12}=1\right]\right|.
\end{array}
\end{equation}
If some coupling (\ref{eq: coupling 2x2}) has this and the analogously
computed maximal values for other equalities in (\ref{eq: equalities}),
then the system is noncontextual: the ``wrong'' influences in it
are all confined to directly changing the distributions of the ``wrong''
random variables. If no such coupling exists, however, the system
is contextual: the influence of $\beta$ upon $A$ and/or $\alpha$
upon $B$ is greater than just distributional changes. As shown in
Dzhafarov, Kujala, and Larsson (2015), Kujala, Dzhafarov, and Larsson
(2015), and Kujala and Dzhafarov (in press), the necessary and sufficient
condition for noncontextuality in accordance with this definition
is 
\begin{equation}
\begin{array}{l}
\max_{k,l\in\left\{ 1,2\right\} }\left|\sum_{i,j\in\left\{ 1,2\right\} }\mathsf{E}\left[A_{i}^{ij}B_{j}^{ij}\right]-2\mathsf{E}\left[A_{k}^{kl}B_{l}^{kl}\right]\right|\\
\leq2+\sum_{i=1}^{2}\left|\mathsf{E}\left[A_{i}^{i1}\right]-\mathsf{E}\left[A_{i}^{i2}\right]\right|+\sum_{j=1}^{2}\left|\mathsf{E}\left[B_{j}^{1j}\right]-\mathsf{E}\left[B_{j}^{2j}\right]\right|.
\end{array}
\end{equation}
For application of this and other criteria of contextuality to available
experimental data in physics and psychology see, respectively, Kujala,
Dzhafarov, and Larsson (2015) and Dzhafarov, Zhang, and Kujala (2015).

\section{Conclusion}

I have argued in this paper that the KPT (Kolmogorovian probability
theory) must allow for stochastically unrelated random variables,
and these must not be confused with stochastically independent ones.
I have argued for radical contextualism: any two random variables
recorded under different conditions (in different contexts) are stochastically
unrelated. There is no fixed set of pairwise stochastically unrelated
random variables: they can be freely introduced and freely coupled.
To couple a given set of stochastically unrelated random variables
means to create their jointly distributed ``copies'' (stochastically
unrelated to the ``originals''). The couplings for a given set of
random variables are typically infinite in number, with no coupling
being ``more correct'' than another. This applies also to couplings
with stochastically independent components. The idea I and Janne Kujala
have been promoting in recent publications is that stochastically
unrelated random variables can be usefully characterized by their
possible couplings, in particular, by determining whether these variables
allow for couplings subject to certain constraints. I have illustrated
this idea on the issue of selective influences, generalized into the
issue of probabilistic contextuality.

\subsection*{Acknowledgments.}

This research has been supported by NSF grant SES-1155956 and AFOSR
grant FA9550-14-1-0318, and A. von Humboldt Foundation. I greatly
benefited from discussions with Matt Jones of the University of Colorado,
my doctoral students Ru Zhang and Victor H. Cervantes, and, of course,
my long-term collaborator (and co-author of the theory presented in
this paper) Janne Kujala. Finally, I should acknowledge my indebtedness
to R. Duncan Luce: my conversations and debates with him have served
as a major source of intellectual inspiration for me for many years.

\section*{References}

\setlength{\parindent}{0cm}\everypar={\hangindent=15pt}Bell,
J. (1964). On the Einstein-Podolsky-Rosen paradox. \emph{Physics}
1, 195-200.

Bohm, D., \& Aharonov, Y. (1957). Discussion of experimental proof
for the paradox of Einstein, Rosen and Podolski. \emph{Physical Review}
108, 1070-1076.

Clauser, J.F., Horne, M.A., Shimony, A., \& Holt, R.A. (1969). Proposed
experiment to test local hidden-variable theories. \emph{Physical
Review Letters} 23, 880-884.

Dzhafarov, E.N. (2003). Selective influence through conditional independence.
\emph{Psychometrika}, 68, 7\textendash 26.,

Dzhafarov, E.N., \& Kujala, J.V. (2010). The Joint Distribution Criterion
and the Distance Tests for selective probabilistic causality. \emph{Frontiers
in Psychology} 1:151, doi: 10.3389/fpsyg.2010.00151.

Dzhafarov, E.N. \& Kujala, J.V. (2012a). Selectivity in probabilistic
causality: Where psychology runs into quantum physics. \emph{Journal
of Mathematical Psychology} 56, 54-63.

Dzhafarov, E.N., \& Kujala, J.V. (2012b). Quantum entanglement and
the issue of selective influences in psychology: An overview. \emph{Lecture
Notes in Computer Science} 7620, 184-195.

Dzhafarov, E.N., \& Kujala, J.V. (2013). Order-distance and other
metric-like functions on jointly distributed random variables. \emph{Proceedings
of the American Mathematical Society} 141, 3291-3301.

Dzhafarov, E.N., \& Kujala, J.V. (2014a). A qualified Kolmogorovian
account of probabilistic contextuality. \emph{Lecture Notes in Computer
Science} 8369, 201-212.

Dzhafarov, E.N., \& Kujala, J.V. (2014b). Contextuality is about identity
of random variables. \emph{Physica Scripta} T163, 014009.

Dzhafarov, E.N., \& Kujala, J.V. (2014c). On selective influences,
marginal selectivity, and Bell/CHSH inequalities. \emph{Topics in
Cognitive Science} 6, 121-128.

Dzhafarov, E.N., \& Kujala, J.V. (2015a). Conversations on contextuality.
In E.N. Dzhafarov, J.S. Jordan, R. Zhang, V.H. Cervantes (Eds). \emph{Contextuality
from Quantum Physics to Psychology}, pp. 1-22. New Jersey: World Scientific
Press.

Dzhafarov, E.N., \& Kujala, J.V. (2015b). Context-content systems
of random variables: The contextuality-by-default theory. arXiv:1511.03516.

Dzhafarov, E.N., \& Kujala, J.V. (in press). Probability, random variables,
and selectivity. In W.Batchelder, H. Colonius, E.N. Dzhafarov, J.
Myung (Eds). \emph{New Handbook of Mathematical Psychology}. Cambridge:
Cambridge University Press.

Dzhafarov, E.N., Kujala, J.V., \& Cervantes, V.H. (2016). Contextuality-by-Default:
A brief overview of ideas, concepts, and terminology. Lecture Notes
in Computer Science 9535, 12-23.

Dzhafarov, E.N., Kujala, J.V., \& Larsson, J.-A. (2015). Contextuality
in three types of quantum-mechanical systems. \emph{Foundations of
Physics} 7, 762-782.

Dzhafarov, E.N., Schweickert, R., \& Sung, K. (2004). Mental architectures
with selectively influenced but stochastically interdependent components.
\emph{Journal of Mathematical Psychology}, 48, 51-64. 

Dzhafarov, E.N., Zhang, R., \& Kujala, J.V. (2015). Is there contextuality
in behavioral and social systems? \emph{Philosophical Transactions
of the Royal Society A} 374, 20150099.

Fine, A. (1982). Hidden variables, joint probability, and the Bell
inequalities. \emph{Physical Review Letters} 48 291-295.

Khrennikov, A. (2005). The principle of supplementarity: A contextual
probabilistic viewpoint to complementarity, the interference of probabilities,
and the incompatibility of variables in quantum mechanics. \emph{Foundations
of Physics} 35, 1655-1693.

Khrennikov, A. (2008). Bell-Boole inequality: Nonlocality or probabilistic
incompatibility of random variables? \emph{Entropy} 10, 19-32.

Khrennikov, A. (2009a). Bell's inequality: Physics meets Probability.
\emph{Information Science} 179, 492-504.

Khrennikov, A. (2009b). \emph{Contextual Approach to Quantum Formalism}.
Berlin: Springer.

Khrennikov, A. (2009c) \emph{Interpretations of Probability}. Berlin:
De Gruyter.

Kolmogoroff, A. N. (1933). \emph{Grundbegriffe der Wahrscheinlichkeitsrech}.
Berlin: Springer Verlag. English translation: Kolmogorov, A.N. (1956).
\emph{Foundations of the Probability Theory}. New York: Chelsea Publishing
Company.

Krantz, D.H., Luce, R.D, Suppes, P., \& Tversky, A. (1971). \emph{Foundations
of Measurement}, vol. I: \emph{Additive and Polynomial Representations}.
New York, NY: Academic Press.

Kujala, J.V., \& Dzhafarov, E.N. (2008). Testing for selectivity in
the dependence of random variables on external factors. \emph{Journal
of Mathematical Psychology}, 52, 128-144.

Kujala, J.V., \& Dzhafarov, E.N. (in press). Proof of a conjecture
on contextuality in cyclic systems with binary variables. \emph{Foundations
of Physics}.

Kujala, J.V., \& Dzhafarov, E.N. (2015). Probabilistic contextuality
in EPR/Bohm-type systems with signaling allowed. In In E.N. Dzhafarov,
J.S. Jordan, R. Zhang, V.H. Cervantes (Eds).\emph{ Contextuality from
Quantum Physics to Psychology}, pp. 287-308. New Jersey: World Scientific.

Kujala, J.V., Dzhafarov, E.N., \& Larsson, J.-A. (2015). Necessary
and sufficient conditions for maximal noncontextuality in a broad
class of quantum mechanical systems. \emph{Physical Review Letters}
115, 150401.

P. Kurzynski, R. Ramanathan, \& D. Kaszlikowski (2012). Entropic test
of quantum contextuality. \emph{Physical Review Letters} 109, 020404.

Larsson, J.-A. (2002). A Kochen-Specker inequality. \emph{Europhysics
Letters}, 58, 799\textendash 805.

Roberts, S., \& Sternberg, S. (1993). The meaning of additive reaction-time
effects: Tests of three alternatives. In D.E. Meyer \& S. Kornblum
(Eds.), \emph{Attention and performance XIV: Synergies in experimental
psychology, artificial intelligence, and cognitive neuroscience} (pp.
611\textendash 654). Cambridge, MA: MIT Press.

Schweickert, R., Giorgini, M., \& Dzhafarov, E.N. (2000). Selective
influence and response time cumulative distribution functions in serial-parallel
networks. \emph{Journal of Mathematical Psychology}, 44, 504\textendash 535.

Schweickert, R., \& Townsend, J. T. (1989). A trichotomy: Interactions
of factors prolonging sequential and concurrent mental processes in
stochastic discrete mental (PERT) networks. Journal of Mathematical
Psychology, 33, 328-347. 

Sternberg, S. (1969). The discovery of processing stages: Extensions
of Donders\textquoteright{} method. In W.G. Koster (Ed.), \emph{Attention
and Performance II. Acta Psychologica}, 30, 276\textendash 315.

Simon, C., Brukner, C., \& Zeilinger, A. (2001). Hidden-variable theorems
for real experiments. \emph{Physical Review Letters,} 86, 4427-4430. 

Svozil, K. (2012). How much contextuality? \emph{Natural Computing,}
11, 261-265.

Thorisson, H. (2000). \emph{Coupling, Stationarity, and Regeneration}.
New York: Springer. 

Townsend, J. T. (1984). Uncovering mental processes with factorial
experiments. \emph{Journal of Mathematical Psychology}, 28, 363-400.

Townsend, J. T., \& Nozawa, G. (1995). Spatio-temporal properties
of elementary perception: An investigation of parallel, serial, and
coactive theories.\emph{ Journal of Mathematical Psychology}, 39,
321-359.

Townsend, J. T., \& Schweickert, R. (1989). Toward the trichotomy
method of reaction times: Laying the foundation of stochastic mental
networks. \emph{Journal of Mathematical Psychology}, 33, 309-327.

von Mises, R. (1957). \emph{Probability, Statistics and Truth}. London:
Macmillan. 

Wang, Z., Busemeyer, J.R. (2013). A quantum question order model supported
by empirical tests of an a priori and precise prediction. \emph{Topics
in Cognitive Science,} 5, 689-710.

Wang, Z., Solloway, T., Shiffrin, R.M., Busemeyer, J.R. (2014). Context
effects produced by question orders reveal quantum nature of human
judgments. \emph{Proceedings of the National Academy of Sciences,}
111, 9431-9436.

Winter, A. (2014). What does an experimental test of quantum contextuality
prove or disprove? \emph{Journal of Physics A: Mathematical and Theoretical,}
47, 42403. 
\end{document}